\documentclass[12pt, a4paper, oneside]{amsart}
\usepackage{amssymb}  
\usepackage{amsmath} 
\usepackage{amsthm}  

\usepackage[cp1251]{inputenc}
\usepackage[english]{babel}
\usepackage{cite, verbatim, array, multirow}
\usepackage[justification=centering]{caption}

\newtheorem{theorem}{Theorem}
\numberwithin{theorem}{section}
\newtheorem{conj}[theorem]{Conjecture}
\newtheorem{cor}[theorem]{Corollary}
\newtheorem{problem}[theorem]{Problem}
\newtheorem{lemma}[theorem]{Lemma}
\theoremstyle{definition}
\newtheorem{example}[theorem]{Example}

\numberwithin{equation}{section}

\newcommand{\Aut}{\operatorname{Aut}}
\newcommand{\Out}{\operatorname{Out}}
\newcommand{\Inn}{\operatorname{Inn}}
\newcommand{\diag}{\operatorname{diag}}
\newcommand{\antidiag}{\operatorname{antidiag}}

\newcommand{\Mod}[1]{\,(\mathrm{mod}\ #1)}

\begin{document}

\title{Finite groups isospectral to simple groups}

\author{\sc Maria A. Grechkoseeva}
\address{Sobolev Institute of Mathematics\\ Koptyuga 4, Novosibirsk 630090, Russia}

\email{grechkoseeva@gmail.com}

\author{Victor D. Mazurov}
\address{Sobolev Institute of Mathematics, Koptyuga 4, Novosibirsk 630090, Russia}
\email{mazurov@math.nsc.ru}

\author{Wujie Shi}
\address{Department of Mathematics, Chongqing University of Arts and Sciences, Chongqing, 402160, China; School of Mathematics, Suzhou University, Suzhou, 215006, China}
\email{shiwujie@outlook.com}

\author{Andrey V. Vasil'ev}
\address{School of Science, Hainan University, Haikou, Hainan, 570228, P.R. China; Sobolev Institute of Mathematics Koptyuga 4, Novosibirsk 630090, Russia}
\email{vasand@math.nsc.ru}

\author{Nanying Yang}
\address{School of Science, Jiangnan University, Wuxi, 214122, China}

\email{yangny@jiangnan.edu.cn}

\thanks{Maria A. Grechkoseeva, Andrey V. Vasil'ev and Nanying Yang were supported by Foreign Experts program in Jiangsu Province (No. JSB2018014). Andrey V. Vasil'ev was supported by the National Natural Science Foundation of China (No. 12171126). Victor D. Mazurov was supported by the RFBR (No. 20-51-00007). Wujie Shi was supported by the the National Natural Science Foundation of China (11171364, 11671063).}

\begin{abstract} The spectrum of a finite group is the set of element orders of this group. The main goal of this paper is to survey results concerning
recognition of finite simple groups by spectrum, in particular, to  list all finite simple groups for which the recognition problem is solved.

\noindent{\bf Keywords:} finite group, simple group, element order, spectrum, recognition by spectrum.
\noindent{\bf MSC:} 20D05, 20D06, 20D60.
 \end{abstract}

\begingroup
\let\MakeUppercase\relax 
\maketitle
\endgroup
\vspace*{-1cm}

\section{Introduction}\label{intro}

The elementary assertion that groups of exponent 2 are abelian and one of the most complicated theorems in the finite group theory on solvability of groups of odd order are similar in the following sense. In both cases, arithmetic, i.e., expressed by numeric parameters, properties of a group allow to conclude on its structure. The strongest conclusion that can be made in this direction is that a group is uniquely, up to isomorphism, determined by a~set of numeric parameters. In such case the group is said to be {\em recognizable} by this set.

The recognition of finite simple groups --- the building blocks of finite group theory --- are definitely of prime interest. According to their classification (CFSG), finite simple groups are exactly
\begin{enumerate}
\item the groups of prime order;
\item the alternating groups of degree at least $5$;
\item the simple classical groups;
\item the simple exceptional groups of Lie type;
\item the 26 sporadic groups.
\end{enumerate}

In 1987, in his letter to J. Thompson, W. Shi conjectured that every finite simple group is recognizable in the class of finite  groups by its order and the set of element orders. In his answer, J. Thompson highly appreciated this conjecture and put forward another one: every finite nonabelian simple group is recognizable in the class of finite groups with trivial center by the set of sizes of conjugacy classes. Later, A. S. Kondrat'ev added the questions about the validity of these conjectures to the {\em Kourovka Notebook} \cite[Problems~12.38, 12.39]{Kou}.

As a result of efforts by numerous mathematicians started in \cite{89Shi1} and \cite{92Chen}, the validity of both conjectures was eventually established: Shi's conjecture --- in 2009 \cite{09VasGrMaz1.t}, and Thompson's conjecture --- in 2019 \cite{19Gor}. Together with \cite[Corollary~5.2]{08KimLuRC}, the validity of Shi's conjecture implies that a finite simple group and an arbitrary finite group having the same Burnside ring are isomorphic, i.e., every finite simple group $G$ is recognizable in the class of finite groups by $G$-sets (Yoshida's problem \cite[Problem 2]{87Yo} for simple groups).

In Shi's conjecture, a recognition is based on the order and the set of element orders. What happens if we leave only the latter of these parameters --- the set of element orders? Following \cite{84Ady.t}, we refer to this set as the {\em spectrum} of a group $G$ and denote it by~$\omega(G)$, while groups with the same spectrum are said to be {\em isospectral}. Observe that another notation for the spectrum is $\pi_e(G)$ (see, for example, \cite{12Shi}).

In the middle of the 1980s W. Shi \cite{84Shi, 86Shi} discovered that $PSL_2(7)$ and $Alt_5$  can be characterized in the class of finite groups solely  by the spectrum, and these results opened a~wide way for investigations of recognizability of groups by spectrum.  Though a finite group $G$ isospectral to a nonabelian simple group $L$ is not necessarily isomorphic to~$L$, it is very close to $L$ in the vast majority of cases. Namely, if $L$ is alternating of sufficiently large degree or classical of sufficiently large dimension, then
$G$ is isomorphic to a group squeezed between $\Inn L$ and $\Aut L$, the groups of inner and of all automorphisms of $L$ respectively. This property of `sufficiently large' simple groups was conjectured by V. D. Mazurov in 2007 and eventually proved in 2015 (see \cite{15VasGr1} and the references therein).

We denote the number of pairwise nonisomorphic finite groups isospectral to a group $G$ by $h(G)$. So $G$ is recognizable (by spectrum) if $h(G)=1$. A group $G$ is said to be {\em almost recognizable} if $h(G)$ is finite, and {\em unrecognizable} if $h(G)=\infty$ \cite{97Shi}. We say that the {\em recognition problem} is solved for $G$ if the number $h(G)$ is known, and if it is finite, then all the groups isospectral to $G$ are described. The main goal of this paper is to survey known results on this problem. The previous surveys can be found in \cite{91Shi, 04Maz, 05MazS, 08GrShiVas, 07Shi, 12Shi, 15Maz} and
partially in the introduction in \cite{15VasGr1}.

Since the study of the recognition problem has naturally gone in parallel with the study of spectra, it is worth noting that the spectra of all finite nonabelian simple groups are known.
Clearly, the spectra of alternating groups are easy to describe, and the spectra of sporadic groups are known due to \cite{85Atlas}. The spectra of linear and unitary groups are found in \cite{08But.t}, and those of symplectic and orthogonal in \cite{10But.t}. For exceptional groups, the task was completed in \cite{18But.t} with describing the spectra of $E_8(q)$.

The structure of the paper is as follows.
Section~\ref{main} is concerned with finite simple groups for which the recognition problem is solved (Theorem \ref{t:main}), and the groups themselves are listed in Appendix. In Section~\ref{general}, we state some general properties of finite groups isospectral to finite simple groups and discuss finite simple groups for which the recognition problem is not solved. Finally, we discuss some related problems and open questions in Section \ref{related}.

\section{Simple groups with solved recognition problem}\label{main}

\subsection{The main theorem}

The finite nonabelian simple groups for which the recognition problem is solved are listed in Tables \ref{tab:l}--\ref{tab:alt} in Appendix. The main result of the section is Theorem \ref{t:main} which
describes finite groups isospectral to $L$ for every simple group $L$ listed in the tables and having $h(L)<\infty$.

We denote the alternating and symmetric groups of degree $n$ by $Alt_n$ and $Sym_n$ respectively. In notation of
the sporadic groups and groups of Lie type we follow
the `Atlas of finite groups' \cite{85Atlas},  with the exception
that we write $^2B_2(q)$, $^2G_2(q)$ and $^2F_4(q)$ for the Suzuki--Ree groups. Also we use the abbreviations $L^\pm_n(q)$ and $E^\pm_6(q)$, where $L_n^+(q)=L_n(q)$, $L_n^-(q)=U_n(q)$, $E^+_6(q)=E_6(q)$, $E_6^-(q)={}^2E_6(q)$. As usual, we identify $L$  with $\Inn L$ and write $\Out L$ for the quotient $\Aut L/\Inn L$, the outer automorphism group of $L$.

\begin{theorem}\label{t:main}
Let $L$ be one of the finite simple groups listed in Tables \ref{tab:l}--\ref{tab:alt}. The following hold:

\begin{enumerate}
\item $h(L)$ is as  in the third column of the tables.

\item  Suppose that  $1<h(L)<\infty$ and $L\neq S_6(2)$, $O_8^+(2)$, $O_7(3)$, $O_8^+(3)$, and let $\Theta$ be a~subset of $\Out L$ specified in the last column of the tables. Then a finite group $G$ is isospectral to $L$ if and only if $L\leq G\leq \Aut L$ and $G/L$ is conjugate in $\Out L$ to $\langle \theta\rangle$ for some $\theta\in \Theta$.

\item Let $\mathcal{L}=\{S_6(2), O_8^+(2)\}$ or $\{O_7(3), O_8^+(3)\}$, and suppose that $L\in\mathcal{L}$.  Then a~finite group $G$ is isospectral to $L$ if and only if $G\in\mathcal{L}$.\end{enumerate}
\end{theorem}

\begin{proof} It should be noted that in this proof, for the sake of brevity, we generally refer not to all relevant papers but only to the papers in which some final steps were made or some remaining cases were settled. We will give a more detailed exposition in Section \ref{general}.

If $h(L)=\infty$, then see the references given in the last column of the tables. If $L$ is an alternating or sporadic, then see \cite{13Gor.t} or \cite{98MazShi} respectively and the references therein.
If $L$ is exceptional, then see \cite[Theorem 3]{16Zve.t}. If $L=L_2(q)$, $L_3(q)$, $U_3(q)$, or $S_4(q)$, then see \cite{94BrShi}, \cite{04Zav1.t}, \cite{06Zav.t}, or \cite{02Maz.t} respectively. If $L$ is a classical
group in characteristic $2$, then see \cite{15VasGr.t}. Finally, if $L=O_7(3)$ or  $O_8^+(3)$, then see \cite{97ShiTan}.

Thus we may assume that $L$ is a classical group in odd characteristic with
$h(L)<\infty$ and $L\neq L_2(q)$, $L_3(q)$, $U_3(q)$, $S_4(q)$, $O_7(3)$,  $O_8^+(3)$. Let $G$ be a finite group isospectral to~$L$. First, we need to prove that $G$ is an almost simple group with socle $L$, that is, $L\leq G\leq \Aut L$.

If $L=L_4^\pm(q)$, then this holds by \cite{20GrZv.t}.
If $L=L_n(q)$, where $n\geqslant 5$ is prime, then see \cite{12GrLyt.t}.
If $L=S_6(q)$, $O_7(q)$, or $O_8^+(q)$, then see \cite{19GrVasZv}. If $L=L_n^\pm(q)$, $S_{2n}(q)$, $O_{2n+1}(q)$, or $O_{2n}^\pm(q)$, and $n\geqslant n_0$ for $n_0$ specified at the top of the corresponding table, then see \cite[Theorem 1.2]{17Sta}.

If $L=S_{2n}(q), O_{2n+1}(q)$, where $n\geqslant 8$ is a~power of 2, or $L=O_{2n}^-(q)$, where $n\geqslant 4$ is a~power of 2, then by \cite[Theorems 1, 2]{09VasGorGr.t}, it follows that $G$ has a~unique nonabelian composition factor and this factor is isomorphic to $L$.
The same is true when $L$ is one of the groups in Table~\ref{tab:dpg} (see \cite{16He}
for $U_n(3)$,  \cite{12HeShi} for $S_{2n}(3)$,  \cite{10SheShiZin.t} for $O_{2n+1}(3)$, and \cite{09HeShi,  09Kon, 08AlKon, 09Kon1, 09AleKon.t} for $O_{2n}^\pm(q)$). Now we apply  \cite[Theorem 1.1]{15Gr} to conclude that the solvable radical of $G$ is trivial, and so $G$ is an almost simple group with socle $L$, as required.

It remains to determine almost simple groups with socle $L$ isospectral to $L$. Such groups are described in \cite{17Gre} for linear or unitary socle, in \cite{16Gr.t} for symplectic or odd-dimensional orthogonal socle, and in
\cite{18Gr.t} for even-dimensional orthogonal socle.\end{proof}

\subsection{Number-theoretical notation and outer automorphisms used in the tables}\label{ss:def}
Given a positive integer $n$, we write $\pi(n)$ and $\tau(n)$ for the set of prime divisors of $n$ and the number of all divisors of $n$ respectively. For a prime $r$, we write $(n)_r$ for the $r$-part of $n$, that is, the highest power of $r$ dividing $n$; and $(n)_{r'}$ denotes the $r'$-part of $n$, that is, the ratio $n/(n)_r$. More generally, if $m$ is a positive integer, then $(n)_m$ is equal to the product of $(n)_r$ over all $r\in\pi(m)$ and $(n)_{m'}=n/(n)_m$. The greatest common divisor and least common multiple of integers $n_1$, $n_2$, \dots, $n_s$ are denoted by $(n_1,n_2,\dots,n_s)$ and $[n_1,n_2,\dots,n_s]$ respectively. For simplicity, we write $(n_1,n_2,\dots,n_s)_r$ instead of $((n_1,n_2,\dots,n_s))_r$ for the $r$-part of $(n_1,n_2,\dots,n_s)$. Also for $\varepsilon\in\{+,-\}$, we write $\varepsilon$
in place of $\varepsilon1$ in arithmetic expressions.

Recall that a Mersenne prime $r$ has the form $r = 2^k - 1$ for some $k$. Following \cite{06Zav.t}, we say that a Mersenne prime $r$ is \emph{special} if $r^2 -r +1$ is also a prime. For example, 3 and 7 are special Mersenne primes (in fact, these are the only \emph{known} special Mersenne primes).

We write $\diag(a_1,a_2,\dots,a_n)$ for the matrix $A = (a_{ij})$ of size $n\times n$ with $a_{ii} = a_i$ for all $i$ and $a_{ij} = 0$ otherwise. We write $\antidiag(a_1,a_2,\dots,a_n)$ for the matrix $A = (a_{ij})$ of size $n\times n$ with $a_{i,n-i+1} = a_i$ for all $i$ and $a_{ij} = 0$ otherwise. The transpose of $A$ is denoted by~$A^\top$.

Let $L$ be a group of Lie type. We are going to introduce notation for outer automorphisms of $L$ involved in the structure of groups isospectral to $L$. Since the Suzuki--Ree groups are recognizable by spectrum, in the following definition of automorphisms, we assume that $L$ is not a~Suzuki--Ree group.

Let $L$ be a group over a field of characteristic $p$ and order $q$.
It is well known that we can identify $L$ with a group of the form
$O^{p'}(\overline L_\sigma)/Z(O^{p'}(\overline L_\sigma)),$
where $\overline L$ is a suitable simple algebraic group over the algebraic closure $\overline F$ of the field of order $p$, $\sigma$ is a suitable
endomorphism of $\overline L$, and $\overline L_\sigma=C_{\overline L}(\sigma)$.
The map $x_\alpha(t)\mapsto x_\alpha(t^p)$, where $x_\alpha(t)$ are root elements of $\overline L$, induces an endomorphism of $\overline L$ denoted by $\varphi_p$ (see \cite[Theorem 1.15.4(a)]{98GorLySol}).
If $q=p^m$, then $\varphi_q$ stands for $(\varphi_p)^m$.

We may assume that $\sigma=\varphi_q\gamma_0$, where $\gamma_0$ is some  graph automorphism. We take $\gamma_0$ to be the graph automorphism induced by a suitable symmetry of the Dynkin diagram of $\overline L$ as in \cite[Theorem 1.15.2(a)]{98GorLySol} with the following two exceptions. If $L$ is $L_n^\varepsilon(q)$ or $O_{2n}^\varepsilon(q)$ with $q$ odd, we need a specific $\gamma_0$ and we define it as follows.

Let $p$ be odd. If $L=L_n^\varepsilon(q)$ with $n\geqslant 3$, then we assume that $\overline L=SL_n(\overline F)$ and define $\gamma_0$ to be the inverse-transpose automorphism $g\mapsto g^{-\top}$ of $\overline L$. If $L=O_{2n}^\varepsilon(q)$,
then we assume that $\overline L=SO(\overline V, f)$, where  $\overline V$ is a vector space of dimension $2n$ over $\overline F$ and $f$ is a nondegenerate quadratic form. We can choose a basis of $\overline V$ so that the matrix of $f$ with respect to this basis is $\antidiag(1,1,\dots,1)$. Denote by $\gamma_0$ the linear transformation of $\overline V$ interchanging the first and last basis vectors and fixing all others.
Now we define $\varphi$ and $\gamma$ to be the images in $\Out L$ of the automorphisms of $L$ induced by $\varphi_p$ and $\gamma_0$ respectively. Also, if $L=L_n^\varepsilon(q)$ and  $\lambda$ is a primitive $(q-\varepsilon)$th root of unity, then $\delta$ is the image in $\Out L$ of the diagonal automorphism of $L$ induced by conjugation by $\diag(\lambda,1,1,\dots,1)$.

If $\alpha\in \Out L$, then $L.\langle \alpha\rangle$ denotes the full preimage of $\langle \alpha\rangle$ in $\Aut L$.

\subsection{Description of the tables and examples of use}

Let $L$ be one of the finite simple groups listed in Tables \ref{tab:l}--\ref{tab:alt}. At the top of each table (except for Table \ref{tab:alt}), we introduce the numbers and automorphisms used in the body of the table.
The first two columns provide restrictions on $L$. The column `$h(L)$'
specifies the number of groups isospectral to $L$. If $h(L)=\infty$, then
the column `Note' provides the reference to the paper where this fact was proved. If $h(L)$ is finite but not 1 and $L\neq S_6(2)$, $O_7(3)$, $O_8^+(2)$,  $O_8^+(3)$, then the column `Note' specifies the subset $\Theta$ of $\Out L$ whose meaning is explained in Theorem~\ref{t:main}. If $L$ is $S_6(2)$, $O_7(3)$, $O_8^+(2)$, or $O_8^+(3)$,
then $h(L)=2$ and the column `Note' specifies the only finite group isospectral but not isomorphic to $L$.
Finally, if $h(L)=1$, then the corresponding entry of this column is empty. Also Tables \ref{tab:l} and \ref{tab:u} contain the references to the following Lemma \ref{l:graph}
which gives explicit arithmetic conditions for some almost simple groups to be not isospectral. These conditions are too voluminous to be conveniently included in the~tables.

\begin{lemma}\label{l:graph}
Let $L=L_n^\varepsilon(q)$, where $n\geqslant 3$, $\varepsilon\in\{+,-\}$ and $q$ is a power of an odd prime~$p$. Then $\omega(L.\langle \gamma\rangle)\neq \omega(L)$ if and only if  one of the following holds:

\begin{enumerate}
 \item  $n=p^{t-1}+2$ with $t\geqslant 1$ and $q\equiv - \varepsilon\pmod{4}$;
 \item  $n=2^{t}+1$ with  $t\geqslant 1$ and $(n,q-\varepsilon)>1$;
 \item $n=p^{t-1}+1$ with  $t\geqslant 1$;
 \item $n$ is even, $(n)_2\leqslant (q-\varepsilon)_2$, and $q\equiv \varepsilon\pmod 4$;
 \item $n$ is even,  $(n)_{2'}>3$, and $(n,q-\varepsilon)_{2'}>1$.%
\end{enumerate}
\end{lemma}

\begin{example} To give an example of use of the tables, let us first determine finite groups isospectral to $L=U_{27}(q)$. If $p=2$ then $h(L)=1$ since $n-1$ is not a power of $2$.

Let $p$ be odd. Again $n-1$ is not a power of $p$, so we apply Lemma \ref{l:graph} to see whether $\omega(L.\langle \gamma\rangle)=\omega(L)$ or not.
It follows that $\omega(L.\langle \gamma\rangle)\neq\omega(L)$ if and only if $p=5$ and $q\equiv 1\pmod 4$. It is clear that the former condition implies the latter.

Let $p=5$. Since $b=((q+1)/(27,q+1), m)_3$ and $(q+1)_3$ is equal to $1$ or  $3(m)_3$ according as $m$ is even or odd, we conclude that $h(L)=1$ if and only if $m$ is even or $(m)_3\leqslant 9$. If $m$ is odd and $(m)_3=3^t\geqslant 27$, then $h(L)=\tau(3^{t-2})=t-1$ and $\omega(G)=\omega(L)$
if and only if $L\leq G\leq L.\langle \psi\rangle$, where $\psi$ is a field automorphism of $L$ of order $3^{t-2}$.

Suppose that $p\neq 2, 5$. Again we need to calculate $b=((q+1)/(27,q+1), m)_3$.  Since $n$ is odd, it follows that $h(L)=\tau(2(m)_2b)\geqslant 2$. Also $\omega(G)=\omega(L)$ if and only if $L\leq G\leq L.\langle \alpha\rangle$, where $\alpha$ is an element of $\langle\varphi\rangle$ of order $2(m)_2b$.
\end{example}

\begin{example} Now we consider the group $L=L_{28}(q)$. If $p=2$, then $h(L)=1$.

Let $p$ be odd.  Since $b=((q-1)/(28,q-1),m)_{28}$ and $(q-1)_2/(4,q-1)\geqslant (m)_2$ if $m$ is even, we see that $(b)_2=(m)_2$.

If $n-1$ is a power of $p$, then $p=3$. Since $n-2\neq 2^s$, we need to look at the value of $(b)_2=(m)_2$. It follows that $h(L)=1$ if $(m)_2\leqslant 2$ and $h(L)=2$ otherwise. In the latter case, $\omega(L)=\omega(L.\langle \chi\gamma\delta\rangle)$, where $\chi$ is a field automorphism of order $2$.

Let $p\neq 2,3$. Suppose that $m$ is odd. Then $b$ is odd, and applying Lemma \ref{l:graph}, we conclude that $\omega(L.\langle \gamma\rangle)\neq\omega(L)$ if and only if $q\equiv 1\pmod 4$ or $q\equiv 1\pmod 7$, in other words, when $(28,q-1)>1$. Thus $h(L)=\tau(b)$ or $2\tau(b)$ depending on whether $(28,q-1)>1$ or not.

Now suppose that $m$ is even. The condition $n=p^s+2^u+1$ holds if and only if $p=5,11,19$ or $23$. In these cases, $h(L)=2\tau(b)-\tau((b)_{2'})$. For example, if $q=5^{42}$, then $b=14$ and $h(L)=6$, with the groups isospectral to $L$ being $L$, $L.\langle \xi\rangle$, $L.\langle \xi\gamma\rangle$, $L.\langle \psi\rangle$, $L.\langle \psi\xi\rangle$, and $L.\langle \psi\xi\gamma\rangle$, where $\psi$ and $\xi$ are field automorphisms of orders $7$ and $2$, respectively. Observe that, contrary to the previous cases,  here we have two maximal nonisomorphic subgroups of $\Aut L$ among the groups isospectral to $L$.

If $p\neq 5,11,19, 23$,  we choose $\kappa\in\{+,-\}$ so that $p\equiv \kappa\pmod 4$. Then $h(L)$ depends on the value of $(p-\kappa)_2$.
For example, if $q=7^2$, then $b=2$ and $h(L)=3\tau(b)-2\tau((b)_{2'})=4,$
and the groups isospectral to $L$ are $L$, $L.\langle \chi\rangle$, $L.\langle \chi\gamma\rangle$, and $L.\langle \chi\delta\rangle$, where $\chi$ is a field automorphism of order $2$.
\end{example}

\section{The structure of groups isospectral to simple groups}\label{general}

\subsection{General results}

We need a definition of the prime graph (or the Gruenberg--Kegel graph) $GK(G)$ of a finite group $G$: the vertex set of this graph is $\pi(G)$ and primes $r,s\in\pi(G)$ are adjacent if and only if $r\neq s$ and $rs\in\omega(G)$. Recall that a coclique is a graph in which all the vertices are nonadjacent. The maximal size of a coclique in $GK(G)$
is denoted by $t(G)$ and for $r\in\pi(G)$, the maximal size of a coclique in $GK(G)$ containing $r$ is denoted by $t(r,G)$.

Observe that the prime graph of a solvable group $G$ contains no cocliques of size 3, and so $t(G)\leqslant 2$. This is a direct consequence of G. Higman's result \cite[Theorem 1]{57Hig}, but as an independent statement, it is due to M. S. Lucido  \cite[Proposition 1]{99Luc}.

The connected components of the prime graphs of finite simple groups were described by J. S. Williams \cite{81Wil} and A. S. Kondrat'ev \cite{89Kon.t}. Explicit criteria of adjacency in these graphs were found by A. V. Vasil'ev and E. P. Vdovin \cite{05VasVd.t} (see also \cite[Section 4]{11VasVd.t} and \cite[Section 2]{12HeShi} for corrections to \cite{05VasVd.t}). One of important corollaries of Vasil'ev and Vdovin's result is that $t(2,L)\geqslant 2$ for all finite simple groups $L$ except some $Alt_n$ with $n\geqslant 27$ \cite[Theorem 7.1]{05VasVd.t}. The importance of this fact for the recognition problem is explained by the following result of A. V. Vasil'ev.

\begin{theorem}[\!\!{\cite{05Vas.t}}]\label{t:structure} Let $G$ be a finite group with $t(2,G)\geqslant 2$. Then $G$ has at most one nonabelian composition factor.
\end{theorem}

In 2013 I. B. Gorshkov \cite{13Gor.t} proved that a finite group $G$ isospectral to $Alt_n$ with  $n\geqslant 5$ and $n\neq 6,10$ is isomorphic to $Alt_n$ and, in particular, has a unique nonabelian composition factor. Also by \cite{04LucMog}, \cite[Proposition 3]{02Maz.t}, \cite[Example 2]{03Ale.t}, \cite{08Sta.t}, and \cite{10Zav.t}, the simple groups isospectral to solvable groups are exactly $L_3(3)$, $U_3(3)$, and $S_4(3)$. Thus we have the following theorem.

\begin{theorem}\label{t:factor}
Let $L$ be a finite nonabelian simple group and suppose that $G$ is a finite group isospectral to $L$. Then $G$ has at most one nonabelian composition factor. Furthermore, a~solvable group isospectral to $L$ exists if and only if $L$ is $L_3(3)$, $U_3(3)$, or $S_4(3)$.
\end{theorem}

Observe that $L_3(3)$, $U_3(3)$ and $S_4(3)$ have disconnected prime graphs, so by the Gruenberg--Kegel theorem \cite[Theorem A]{81Wil}, a solvable group
isospectral to one of them is Frobenius or 2-Frobenius.

In what follows in this section we assume that $G$ is a nonsolvable group isospectral to a simple group. By Theorem \ref{t:factor}, it follows that $G$ has the normal series \begin{equation} \label{e:str} 1\leq K<H\leq G,\end{equation} where $K$ is the solvable radical of $G$, that is, the largest normal solvable subgroup of $G$, $S=H/K$ is a nonabelian simple group, and $G/K$ is a subgroup of $\Aut S$. By the validity of Mazurov's
conjecture, for `sufficiently large' simple groups, we have $S\simeq L$ and $K=1$. But what is known about $S$, $K$, and $G/H$, in general?

If $L$ has a disconnected prime graph, then it follows from the Gruenberg--Kegel theorem \cite[Theorem A]{81Wil} and Thompson's theorem on nilpotency of Frobenius kernel that $K$ is nilpotent. It turns out that $K$ is nilpotent for all $L$ other than $Alt_{10}$.

\begin{theorem}[\!\!{\cite{20YanGrVas}}]\label{t:nilpotent}
Let $L$ be a finite nonabelian simple group and suppose that $G$ is a~nonsolvable finite group isospectral to $L$. If $L\neq Alt_{10}$, then the solvable radical of $G$ is nilpotent.
\end{theorem}

By \cite{98Maz.t}, the group $Alt_{10}$ is isospectral to a group of the form $(7^4 \times 3^{12}):(2.Sym_5)$ whose solvable radical is a Frobenius group with complement of order $2$, so $Alt_{10}$ is indeed an exception to Theorem \ref{t:nilpotent}.
Another exceptional feature of this example is that the solvable radical here is a $3$-primary group. Also by \cite{19Lyt}, there is a nonsolvable group with biprimary solvable radical isospectral to $S_4(3)$. We suggest that these are the only examples of a non-primary solvable radical.

\begin{conj} \label{c:r-group}
Let $L$ be a finite nonabelian simple group and suppose that $G$ is a nonsolvable finite group isospectral to $L$. Suppose that $L\neq S_4(3), Alt_{10}$. If $K$ is the solvable radical of $G$, then $|\pi(K)|\leqslant 1$.
\end{conj}

Recall that $G/H$ can be identified with a subgroup of $\Out S$, and so if $S$ is sporadic, or $Alt_n$ with $n\neq 6$, or a Suzuki--Ree group, then $G/H$ is cyclic.
Recently, it was proved that $G/H$ is always cyclic.

\begin{theorem}[\!\!{\cite[Theorem 1]{21GreVas_arxiv}}]\label{t:cyclic}
Let $L$ be a finite nonabelian simple group and suppose that $G$ is a nonsolvable finite group isospectral to $L$. If $1\leq K<H\leq G$ is the series as in \eqref{e:str}, then $G/H$ is cyclic.
\end{theorem}

\subsection{The groups $K$ and $G/H$ in the case $S\simeq L$}
The structure of  $K$ and $G/H$ is best studied  when $S\simeq L$. We refer to a finite group that can be homomorphically mapped onto $L$ as a cover of $L$. If $S\simeq L$, then not only $G$ but also both $H$ and $G/K$ in \eqref{e:str} are isospectral to $L$. So we have, on the one hand, a cover of $L$ isospectral to $L$ and, on the other hand, an almost simple group with socle $L$ isospectral to $L$.

We say that $L$ is recognizable by spectrum among covers (automorphic extensions) if every cover of $L$ (almost simple group with socle $L$) isospectral to $L$ is isomorphic to $L$. The alternating groups are recognizable among covers by \cite{99ZavMaz.t} and among automorphic extensions by an easy number-theoretical argument. The following groups are also recognizable among covers and automorphic extensions: the sporadic groups \cite{87Shi, 87Shi2, 93ShiLi, 98MazShi, 94PrShi, 88Shi, 89Shi, 90Shi, 89ShiLi, 94Shi}, the Suzuki--Ree groups \cite{92Shi, 93BrShi, 99DenShi, 00LipShi}, $L_2(q)$ \cite{87Shi2, 94BrShi}, $G_2(q)$ \cite{13VasSt.t}, $F_4(2^m)$ \cite{04CaoChenGr.t}, $E_7(q)$\cite{14VasSt.t}, and $E_8(q)$ \cite{10Kon}.

The first examples of groups not recognizable among covers, namely $U_3(3)$ and $U_3(7)$, were found by V. D. Mazurov in \cite{98Maz.t}. A. V. Zavarnitsine \cite{06Zav.t} generalized these examples to the assertion  that $U_3(q)$ is recognizable among covers if and only if  $q$ is not a special Mersenne prime (recall from Section \ref{main} that a Mersenne prime $r$ is special if $r^2-r+1$ is a prime too).
The covers of other unitary and linear groups were settled in \cite{04Zav1.t,07ZavMaz, 08Zav1.t, 11Gr, 20GrSk.t}, of symplectic and orthogonal groups in \cite{11Gr, 15Gr}. The groups $^3D_4(q)$ with $q\neq 2$, $F_4(q)$,  $E_6^\varepsilon(q)$ and $E_7(q)$ are recognizable among covers by \cite{15Gr, 14VasSt.t}. Finally, $^3D_4(2)$ is not recognizable among covers by \cite{13Maz.t}. The following theorem summarizes the results on recognition among covers.

\begin{theorem}\label{t:cover}
Let $L$ be a finite nonabelian simple group and $G$ is a finite group such that $G/K\simeq L$ for some nontrivial normal subgroup $K$ of $G$. If $L\neq U_3(q)$, where $q$ is a special Mersenne prime, and $L\neq U_5(2)$, $^3D_4(2)$, then $\omega(G)\neq \omega(L)$.
\end{theorem}

`Atlas of finite groups' provides several examples of almost simple groups isospectral to their socles, for example, $L_3(5)$ is isospectral to its extension by the graph automorphism.
It turned out that this situation is rather typical for simple groups of Lie type, in other words, simple groups of Lie type are mostly not recognizable among automorphic extensions. This leads to the problem of describing almost simple groups of Lie type isospectral to its socle. A. V.  Zavarnitsine \cite{04Zav,06Zav.t} solved this problem for $L_3(q)$ and $U_3(q)$ and provided the first examples showing that there can be arbitrarily many almost simple groups with the same socle isospectral to this socle. Later the problem was solved for the classical groups in characteristic 2  \cite{08Gr.t, 13GrShi, 14Zve} and  the groups $^3D_4(q)$, $F_4(q)$, $E_6^\varepsilon(q)$, $E_7(q)$ \cite{16GrZv, 16Zve.t}. The final steps relating to classical groups in odd characteristic were done by M. A. Grechkoseeva \cite{16Gr.t, 17Gre, 18Gr.t}, with substantial use of the description of the spectra of these groups given by A. A. Buturlakin \cite{08But.t, 10But.t}.
As a result of all these contributions, we have the following theorem.

\begin{theorem}\label{t:aut}
Let $L$ be a finite nonabelian simple group. The groups $G$ such that $L\leq G\leq \Aut L$ and $\omega(G)=\omega(L)$ are known.
\end{theorem}

\subsection{Quasirecognizability}
Discussing the factor $S=H/K$, it is convenient to use the notion of
quasirecognizability introduced by Kondrat'ev \cite{02AlKon.t}. A finite simple group $L$ is said to be \emph{quasirecognizable by spectrum} if every finite group $G$ isospectral to $L$ has a unique nonabelian composition factor and this factor is isomorphic to $L$. So in terms of \eqref{e:str}, $L$ is not quasirecognizable if there exists some $G$ isospectral to $L$ in which $S\not\simeq L$.

All sporadic, alternating and exceptional groups of Lie type, other than $A_6$, $A_{10}$ and $J_2$, are quasirecognizable by spectrum. This was established in  \cite{87Shi, 87Shi2, 93ShiLi, 98MazShi, 94PrShi, 88Shi, 89Shi, 90Shi, 89ShiLi, 94Shi} for sporadic groups,  \cite{00KonMaz.t, 00Zav.t, 10ShaJia, 13Gor.t} for alternating groups, \cite{92Shi, 93BrShi, 99DenShi} for Suzuki--Ree groups, \cite{13VasSt.t, 03AlKon.t} for $G_2(q)$, \cite{05AlKon.t, 06Al.t} for $^3D_4(q)$ (see also the description of groups isospectral to $^3D_4(2)$ in  \cite{13Maz.t}), \cite{03AlKon.t, 05AlKon.t} for $F_4(q)$, \cite{07Kon.t} for $E_6^\pm(q)$, \cite{14VasSt.t} for $E_7(q)$, and \cite{02AlKon.t} for $E_8(q)$.

The following classical groups over fields of characteristic 2
are not quasirecognizable by spectrum (in the brackets we provide an example of $S\not\simeq L$): $U_4(2)$ ($S=Alt_5$ \cite{98Maz.t}), $U_5(2)$ ($S=M_{11}$ \cite{98Maz.t}), $S_4(2^m)$ ($S=L_2(2^{2m})$ \cite{00MazXuCao.t}), $S_8(2^m)$ ($S=O_8^-(2^m)$ \cite{15VasGr1}), $S_6(2)$ and $O_8^+(2)$ (these groups are isospectral). All the other classical groups over fields of characteristic 2 are quasirecognizable (see  \cite{87Shi, 00MazXuCao.t, 08MazChe.t, 08VasGr.t, 10Gr, 10Xu} for linear and unitary groups and \cite{12Sta.t, 15VasGr.t, 15VasGr1, 15Vas} for symplectic and orthogonal groups).

Suppose that $q$ is odd. The groups $L_2(q)$, $L_3(q)$ and $U_3(q)$, except for $L_2(9)$, $L_3(3)$ and $U_3(3)$, are quasirecognizable. The groups $S_4(q)$ are quasirecognizable if and only if $q=3^{2k+1}>3$. For the classical groups in larger dimensions, the study of $S$ fell naturally into four cases according as $S$ is sporadic, or alternating, or of Lie type in characteristic $p$, or of Lie type in characteristic not $p$. The first three cases are covered in \cite{09VasGrMaz.t, 11VasGrSt.t, 12Sta.t, 15VasGr1}. The final result concerning these three cases is the following.

\begin{theorem}\label{t:factorS}
Let $L$ be a finite simple classical group other than $L_2(9)\simeq Alt_6$, $L_4(2)\simeq Alt_8$, $L_3(3)$, $U_3(3)$, $U_4(2)\simeq S_4(3)$, $U_5(2)$. Suppose that $G$ is a finite group isospectral to $L$ and $S$ is the unique nonabelian composition factor of $G$.
Then the following holds:
\begin{enumerate}
 \item $S$ is neither alternating nor sporadic;
 \item if $S$ is a group of Lie type in the same characteristic as $L$ then
 either $S\simeq L$ or one of the following holds:
 \begin{enumerate}
  \item $L=S_4(q)$, where $q\neq 3^{2k+1}$, and $S=L_2(q^2)$;
  \item $L=S_8(q)$ or $O_9(q)$, and $S=O_8^-(q)$;
  \item $\{L,G\}=\{S_6(2), O_8^+(2)\}$;
  \item $\{L,G\}=\{O_7(3), O_8^+(3)\}$.
 \end{enumerate}
\end{enumerate}
\end{theorem}

Items (a), (c) and (d) of (ii) in Theorem \ref{t:factorS} indeed give rise to exceptions, that is, there is a group $G$ isospectral to the corresponding $L$ and having the corresponding $S$ as a composition factor. The same is true for Item (b) if we exclude $L=S_8(7^m)$, for which the question is still open.

\begin{problem}\label{p:s8} Suppose that $L=S_8(7^m)$ and $S=O_8^-(7^m)$. Determine for which $m$ there exists a finite group $G$ isospectral to $L$ and having $S$ as a composition factor.
\end{problem}

In fact, the recognition problem for $S_{8}(7^m)$ is equivalent to Problem \ref{p:s8}. Indeed,
suppose that we fix $m$ and let $L=S_8(7^m)$. By \cite[Lemma 4.3]{18Gr.t}, if $G$ is isospectral to $L$ then $G$ is either an almost simple group with socle $L$, in which case applying Theorem~\ref{t:aut} we see that $G=L$, or an extension of a nontrivial $7$-group by an almost simple group with socle $O_8^-(7^m)$. So if the group $G$ of Problem \ref{p:s8} does not exist, then  $h(L)=1$, and if it exists, then
$h(L)=\infty$ by Theorem \ref{t:shimaz} below.

The last case when $S$ is a group of Lie type in characteristic not $p$ is the most difficult, and
for a long time this case could be excluded only in some special circumstances, for example, when the prime graph of $L$ is disconnected. In 2015 A. V. Vasil'ev \cite{15Vas} made a~breakthrough in this direction eliminating this case for all classical groups of  dimension at least 62 (more precisely, for all classical groups $L$ with $t(L)\geqslant 23$). Later A. S. Staroletov \cite{17Sta} lowered the bound  to 38. More precisely, Staroletov's bound depends on the type of $L$ and is equal to $27$ for $L=L_n^\pm(q)$, $32$ for $L=S_{2n}(q), O_{2n+1}(q)$, 38 for $L=O_{2n}^+(q)$ and 36 for $O_{2n}^-(q)$ (cf. the number $n_0$ given in Tables \ref{tab:l}--\ref{tab:o-}).

\subsection{Simple groups for which the recognition problem is not solved}
Suppose that $L$ is a  simple group for which the recognition problem is not solved, $G$ is a finite group isospectral to $L$ and $S$ is the unique nonabelian composition factor of $G$.
By Theorem \ref{t:main} and the preceding part of this section, $L$ is a classical group in odd characteristic
and the open part of recognition problem is quasirecognizability of $L$. More precisely, either
$L=S_8(7^m)$ and we refer the reader to the remark after Problem \ref{p:s8}, or $L\neq S_8(7^m)$ and the question is whether $S$ can be a group of Lie type in characteristic different from the characteristic of $L$. We conjecture that the answer to this question is negative.

\begin{conj}\label{c:rest} Let $q$ be odd and $L$ be one of the following simple groups:
\begin{enumerate}
\item $L_n(q)$, where $6\leqslant n\leqslant 26$ and $n$ is not prime;
\item $U_n(q)$, where $5\leqslant n\leqslant 26$;
\item $O_{2n}^+(q)$, where $5\leqslant n\leqslant 18$;
\item $O_{2n}^-(q)$, where $5\leqslant n\leqslant 17$ and $n\neq 8,16$;
\item $S_{2n}(q)$ and $O_{2n+1}(q)$, where $5\leqslant n\leqslant 15$ and $n\neq 8$.
\end{enumerate}
If $G$ is a finite group isospectral to $L$, then the unique nonabelian composition
factor $S$ of $G$ cannot be a group of Lie type in characteristic coprime to $q$.
\end{conj}

There are some partial results on Conjecture \ref{c:rest}:  $S$ is not $E_8(u)$, $E_7(u)$ by \cite{21Sta.t} and $S$ is not $L_2(u)$, $S_4(u)$, $G_2(u)$, $^2B_2(u)$, $^3D_4(u)$ by \cite[Proposition 4.8]{20GrZv.t}. Also by Table \ref{tab:dpg}, the recognition problem is solved for some of the groups $L$ in question with $q=3$ or $5$.

Observe that in the proof of Theorem \ref{t:main} for $L=L_n^\pm(q)$, $S_{2n}(q)$, $O_{2n+1}(q)$, or $O_{2n}^\pm(q)$, the condition $n\geqslant n_0$ is necessary to prove that $S\simeq L$ but it does not influence the description of almost simple groups with socle $L$ isospectral to $L$ by \cite{16Gr.t, 17Gre, 18Gr.t}. Thus Conjecture \ref{c:rest} is in fact is equivalent to

\begin{conj} \label{c:final} Theorem {\rm \ref{t:main}} remains valid if we replace $n_0$ by $5$ in Tables {\rm \ref{tab:l}--\ref{tab:o-}}.
\end{conj}

Conjecture \ref{c:final} covers all nonabelian simple groups except $S_8(7^m)$, so the validation of this conjecture and the solution of Problem \ref{p:s8} will finalize the study of recognition of simple groups by spectrum.

\subsection{Unrecognizable simple groups} The necessary and
sufficient condition for a finite group to be unrecognizable by spectrum was established by V. D. Mazurov and W. Shi \cite{12MazShi.t}.

\begin{theorem}[\!\!{\cite{12MazShi.t}}]\label{t:shimaz}
Suppose that $G$ is a finite group.
\begin{enumerate}
\item If $V$ is an elementary abelian normal subgroup of $G$ and $G_1=V\rtimes G$ is the natural
semidirect product of $V$ and $G$ with $G$ acting by conjugation, then $\omega(G_1)=\omega(G)$. In particular, if $G$ has a nontrivial solvable normal subgroup, then $h(G)=\infty$.
\item If $h(G)=\infty$, then $G$ is isospectral
to a finite group with nontrivial normal solvable subgroup.
\end{enumerate}
\end{theorem}

Suppose that a simple group $L$ is isospectral to a finite group $G$ with nontrivial solvable radical. Using Theorem \ref{t:shimaz}, we can construct infinitely many finite groups isospectral to $L$ but it is not quite clear how to describe all finite groups isospectral to $L$.

One way is to consider not all the groups isospectral to $L$ but only minimal in some sense. In the same work \cite{12MazShi.t},
V. D. Mazurov and W. Shi introduced the following notion: if $\omega$ is  a finite set of natural numbers, then a finite group $G$ is said to be \emph{$\omega$-critical} if $\omega=\omega(G)$ and $\omega\neq\omega(H)$ for any proper section $H$ of $G$ (that is, for any $H=A/B$, where $1\leq A\trianglelefteq B\leq G$ and either $A\neq 1$, or $B\neq G$). Also they proved that for every $\omega$, there are only finitely many $\omega$-critical finite groups and at most one simple $\omega$-critical group.
Given an unrecognizable simple group $L$, we refer to $\omega(L)$-critical groups as just critical. The critical groups are known for $Alt_6$ \cite{13LytY}, $Alt_{10}$ \cite{15LytY.t}, $J_2$ \cite{15LytY.t}, $^3D_4(2)$ \cite{13Maz.t} and $L_3(3)$ \cite{13LytY}. Some critical groups are known for $U_3(3)$ \cite{17Lyt} and $S_4(3)$ \cite{19Lyt}.

We choose to describe groups isospectral to unrecognizable simple groups in the following way. In the first column of Table \ref{tab:un}, we list all known unrecognizable simple groups~$L$. It is clear that there are almost simple groups with socle $L$ isospectral to $L$, and these are known by Theorem \ref{t:aut}, so in the next three columns we give information about groups $G$ that are isospectral to $L$ but are not almost simple groups with socle $L$. For solvable groups $G$, we just indicate whether $G$ is Frobenius or 2-Frobenius.  For nonsolvable $G$, we list all possible combinations of $\pi(K)$, $S$ and $G/H$. In the column `Examples', for every combination, we provide an example of such a group $G$ in notation of \cite{85Atlas}.

It should be noted that the ristrictions on $\pi(K)$, $S$ and $G/H$ in the case $L=O_9(2^m)$ are known but there are no published proofs of them. Also we note that the list of unrecognizable groups in \cite[Table 1]{15Maz}, as well as that in \cite[Table 1]{17Lyt}, includes the group $L_4(13^{24})$ due to \cite{08Zav2}. But, in fact, this group is almost recognizable by \cite{20GrSk.t, 20GrZv.t}.

Finally, we remark that under Conjecture \ref{c:rest}, every unrecognizable simple group is either listed in Table \ref{tab:un}, or isomorphic to $S_8(7^m)$ (cf. Problem \ref{p:s8}).

\section{Related questions}\label{related}

\subsection{Recognition outside the class of finite groups}
There are some natural questions related to the  problem of recognizing simple groups by spectrum. The first question arises when we enlarge the class of groups in which we seek groups isospectral to a given simple group.

A natural extension of the class of finite groups is the class of locally finite groups. Theorem \ref{t:shimaz} has
a simple but useful corollary concerning recognizability of groups in the class of locally finite groups.

\begin{cor}\label{cor:loc}
Let $G$ be a finite group.
\begin{enumerate}
 \item If $h(G)<\infty$, then every locally finite group isospectral to $G$ is finite.
 \item If $h(G)=\infty$, then there is an infinite locally finite group isospectral to $G$.
\end{enumerate}
\end{cor}

\begin{proof}
 (i) Suppose that there is a locally finite but not finite group $H$ isospectral to $G$.
For every $a\in \omega(G)$, let $g_a\in H$ be of order $a$. Then $H_0=\langle g_a\mid a\in\omega(G)\rangle$ is a finite group isospectral to $G$. Taking $H_i=\langle H_{i-1}, x_i\rangle$ with $x_i\in H\setminus H_{i-1}$  for $i=1,2,\dots$,
we obtain infinitely many finite groups $H_i$ isospectral to $G$, contrary to the assumption $h(G)<\infty$.

(ii) By Theorem \ref{t:shimaz}, there is a finite group $H$ isospectral to $G$ and having
a nontrivial  elementary abelian normal subgroup $V$. Furthermore, we have an ascending
chain of finite groups $H_i$ isospectral to $G$:
$$H_0=H\leq H_1=V\rtimes H_0\leq H_2=V\rtimes H_1\leq \dots.$$
The group $\cup_{i=1}^\infty H_i$ is an infinite locally finite group isospectral to $G$.
\end{proof}

By Corollary \ref{cor:loc}, if a finite group $G$ is (almost) recognizable in the class of finite groups then $G$ is (almost) recognizable in the class of locally finite groups. This becomes false for the class of periodic groups. For example, by \cite[Theorem 1]{15MOS.t}, there exist infinitely many prime powers $q$ such that $L_2(q)$ is isospectral to an infinite group (recall from Section~\ref{main}, that all $L_2(q)$, except for $L_2(9)$, are recognizable in the class of finite groups).
The deep and difficult problem of recovering a periodic group from its spectrum (it is worth remembering that the finiteness of the Burnside group $B(2,5)$ is still in question) lies outside the scope of this survey, and we refer the reader to \cite{09MazShi, 14LytMaz, 18HLM}.

\subsection{Recognition of almost simple groups} The second question is what non-simple finite groups are recognizable or almost recognizable by spectrum. Recall that the socle $Soc(G)$
of a finite group $G$ is the  subgroup generated by all minimal normal subgroups of~$G$. By Theorem \ref{t:shimaz}, if $G$ is an almost recognizable finite group, then
\begin{equation}\label{eq:structure}
L_1\times L_2\times \dots\times L_k=Soc(G)\leq G\leq \Aut (L_1\times L_2\times \dots\times L_k)
\end{equation}
for some nonabelian simple groups $L_i$, $1\leqslant i\leqslant k$.
Assume first that $k=1$ in \eqref{eq:structure}, that is, $G$ is an almost simple group.

The recognition problem is solved for all symmetric groups $Sym_n$ except $Sym_{10}$ (see \cite{14Gor.t, 16GorGri.t} and references therein). Namely, they are unrecognizable if $n\in\{2,3,4,5,6,8\}$ and recognizable otherwise.

\begin{problem}
Is the symmetric group $Sym_{10}$ recognizable by spectrum?
\end{problem}

By \cite{04MogZokDar, 15MazMog}, if $L$ is a sporadic group and $L\neq J_2$, then $\Aut L$ is recognizable. Observe that $h(\Aut J_2)$ is either $1$ or $\infty$ but the exact value of this number is still not known  (see \cite{17ZhuShe} for more details).

\begin{problem}
Is the automorphism group of the sporadic group $J_2$ recognizable by spectrum?
\end{problem}

We now turn to almost simple groups of Lie type. Clearly, numerous examples of almost recognizable almost simple groups of this sort are given by Theorem \ref{t:main}, namely, these are non-simple almost simple groups isospectral to simple groups $L$ of Item (ii) of this theorem. On the other hand, applying known facts on spectra of almost simple groups, one can easily find infinitely many unrecognizable almost simple groups. For example, let $q$ be a power of a prime $p$, $n$ an odd number such that $n-1$ is not a power of $p$ and $G=PGL^\varepsilon_n(q)$. Then $\omega(G)=\omega(SL^\varepsilon_n(q))$ \cite[Proposition~2.6]{08But.t}. If, in addition, $(n,q-\varepsilon1)\neq1$, then the center of $SL^\varepsilon_n(q)$ is nontrivial, so the group $G$ is unrecognizable \cite[Corollary~4]{08But.t}. One  of the natural obstacles to solving the recognition problem for almost simple groups is that despite some progress (see, for example,  \cite{17Gr.t}), the following problem is still open.

\begin{problem}
Give an explicit description of the spectra of finite almost simple groups of Lie type.
\end{problem}

Even when $Soc(G)=L_2(q)$, the only general result is that $h(PGL_2(q))=\infty$ if $q$ is a~prime or $q=9$, and $h(PGL_2(q))=1$ otherwise \cite{07CheMazShi}. Also the problem is solved for $q=4,5$, in which case $G= Sym_5$, and for $q=9$, in which case  $G$ is one of the groups $Sym_6$, $PGL_2(9)$, $M_{10}$, and $\Aut (L_2(9))$. By the above, $h(Sym_6)=h(PGL_2(9))=\infty$. By \cite{91BrShi}, we have $h(M_{10})=1$. Since $\omega(\Aut (L_2(9)))=\omega(SL_2(9))$, it follows that $h(\Aut (L_2(9)))=\infty$.

\begin{problem}\label{p:l2(q)}
Solve the recognition problem for all almost simple groups with socle $L_2(q)$.
\end{problem}

\subsection{Recognition of groups with non-simple socle}

Let now $G$ be a group satisfying~\eqref{eq:structure} with $k>1$. Generally, one can hardly expect the group $G$ to be recognizable or almost recognizable. For example, let $G_i=Alt_5\times Alt_{4+i}$ with $i=1,2$. Clearly, $\mu(G_i)=\{2^i\cdot3, 2^i\cdot5, 3\cdot5\}$,
where $\mu(G)$ is the subset of $\omega(G)$ consisting of numbers maximal under divisibility.
So $G_i$ is unrecognizable being isospectral to the direct product of the Frobenius groups of the form  $5^2:2^i$ and $5^2:3$ (their kernels are elementary abelian groups of order 25 and complements are cyclic group of orders~$2^i$ and~$3$ respectively). Nevertheless, if we do not presuppose that the simple factors of the socle of $G$ are isomorphic to each other, then we can construct a recognizable group with $k$ arbitrarily large.

Define primes $p_i$ and groups $L_i$ in the following way.
\begin{equation}\label{eq:example}
 \text{\begin{minipage}{0.92\textwidth}Let $p_1=7$ and $L_1={}^2B_2(2^{p_1})$. For $i>1$, define $p_{i}$ to be the smallest prime larger than $p_{i-1}$ and not lying in $\cup_{j<i}\pi(L_j)$, and  put $L_{i}={}^2B_2(2^{p_i})$\end{minipage}}
\end{equation}

The order, the spectrum and subgroups of  the simple Suzuki groups are well-known and can be found in \cite{62Suz}. We have $|{}^2B_2(q)|=q(q-1)(q^{2}+1)$, and so $\pi(L_i)\cap\pi(L_j)=\{2,5\}$ if $i\neq j$. Also $p_i\not\in\pi(L_j)$ for all $i$ and $j$.
Since $\mu({}^2B_2(q))=\{4, q-1, q-\sqrt{2q}+1, q+\sqrt{2q}+1\},$
if we take $r_{i1}\in\pi(q_i-1)$, $r_{i2}\in\pi(q_i-\sqrt{2q_i}+1)$ and $r_{i3}\in\pi(q_i+\sqrt{2q_i}+1)$, where $q_i=2^{p_i}$, then $\{r_{i1}, r_{i2}, r_{i3}\}$ is a coclique in $GK(L_i)$. Observe that we can take $r_{i2}$ and $r_{i3}$ to be larger than $5$ because $p_i\geqslant 7$.

Now let $G_k=\prod_{i=1}^k{}^2B_2(q_i)$ and let $G$ be a finite group isospectral to $G_k$.
It is clear that $\{r_{i1}, r_{i2}, r_{i3}\}$ is still a coclique in $GK(G_k)$, so $t(G)\geqslant 3$, and hence $G$ is not solvable by \cite{57Hig}. Let $K$ be the solvable radical of $G$
and let $S_1\times \dots \times S_t$ be the socle of $G/K$. Since $3\not\in\pi(G_k)$, we see that every $S_i$ is equal to ${}^2B_2(2^{m_i})$ for some $m_i$. By the Bang--Zsigmondy theorem \cite{86Bang}, for every $m\geqslant 7$, there is a prime dividing $2^m-1$ but not $2^l-1$ for $l<m$. Using this fact,
it is not hard to see that $\{m_1,\dots,m_t\}\subseteq\{p_1,\dots,p_k\}$. If $m_i=m_j=p_l$ for some $i\neq j$, then $r_{l1}r_{l2}\in\omega(G)\setminus\omega(G_k)$, so $m_i\neq m_j$ if $i\neq j$.
In particular, $G/K\leqslant \Aut S_1\times \dots\times \Aut S_t$. Observe that $|\Out L_i|=p_i$, so $G/K=Soc(G/K)$. If some $p_i\not\in\{m_1,\dots,m_t\}$, then $\{r_{i1}, r_{i2}, r_{i3}\}\subseteq\pi(K)$ and $t(K)\geqslant 3$, which is a contradiction since $K$ is solvable. Thus $G/K$ is isomorphic to $G_k$, and we may assume that $S_i\simeq L_i$ for all $1\leqslant i\leqslant k$.

Suppose that $K\neq 1$. There are a prime $r$ and a normal subgroup $K_0$ of $K$ such that $V=K/K_0$ is an elementary abelian $r$-group. For every $1\leqslant i\leqslant k$, let $x_i$ be an element of $G$ of order $r_{i2}$ if $r$ not dividing $q_i-\sqrt{2q_i}+1$, or of order $r_{i3}$ otherwise.  Let $V_0=V$ and $V_1=C_V(x_1)$. Then $S_2$ acts on $K_1$ and we define $V_2=C_{V_1}(x_2)$ and so on. If $V_k\neq 1$, then $G$ contains an element of prohibited order $r|x_1|\dots|x_k|$. Suppose that $V_k=1$ and let $i$ be the smallest $i$ such that $V_i=1$. It is clear that the image of $x_i$ in $G/K$ lies in $S_i$, so by
\cite[Lemma 3.8]{20GrZv.t}, it follows that $r=2$. The  group $S_i$ includes a Frobenius subgroup with kernel $\langle x_i\rangle$ and cyclic complement of order $4$. Since $C_{V_{i-1}}(x_i)=1$, by \cite[Lemma 1]{97Maz.t}, we conclude that $G$ has an element of order $8$, which is a contradiction. Thus $K=1$, and so we prove the following

\begin{theorem} For every $k\geqslant 1$, the group $\prod_{i=1}^k{}^2B_2(2^{p_i})$, where $p_i$ are the primes as in \eqref{eq:example}, is recognizable by spectrum.
\end{theorem}

Let now $L=L_1\simeq\ldots\simeq L_k$ and so $Soc(G)=L^k$.  We begin with the case when $G=Soc(G)$, that is, $G$ is a product of nonabelian simple groups. The following theorem collects all known examples of recognizable groups of this type.

\begin{theorem} The following groups are recognizable by spectrum:
\begin{enumerate}
 \item ${}^2B_2(2^7)\times {}^2B_2(2^7)$ {\rm \cite{97Maz1.t}};
  \item $J_4\times J_4$ {\rm \cite{21GorMas}};
 \item $L_{2^m}(2)\times L_{2^m}(2)\times L_{2^m}(2)$ with $m\geqslant 6${\rm \cite{20Gor}}.
\end{enumerate}
\end{theorem}

Note that for every finite simple group $L$, there exists $k$ (depending on $L$) such that $G=L^k$ is unrecognizable. Indeed, if $k\geq|\mu(L)|$, then $\mu(G)=\{m\}$ is a one-element set and so $H$ is unrecognizable being isospectral to a cyclic group of order $m$. For example, $\mu(Alt_5)=\{2,3,5\}$, so $(Alt_5)^3$ is isospectral to a cyclic group of order~$30$ (as we remarked above, $(Alt_5)^2$ is also unrecognizable). The question about existence of a global constant $c$ such that $L^c$ is unrecognizable for all finite simple groups $L$ is open. We put it as follows.

\begin{problem}\label{p:nonsimple}
Does there exist a finite recognizable by spectrum group which is a direct product of $k$ copies of a finite nonabelian simple group for arbitrarily large $k$?
\end{problem}

The case $G\neq Soc(G)$ is even more mysterious. In \cite{94Shi} W. Shi conjectured that the group $G$ having a minimal normal subgroup $N$ with $|\mu(N)|=1$ must be unrecognizable. However, as shown by V. D. Mazurov in \cite{97Maz1.t}, if $G$ is the permutation wreath product of the simple Suzuki group $B={}^2B_2(2^7)$ and a subgroup of the symmetric group on $23$ letters which is isomorphic to a Frobenius  group  of order $23\cdot 11$, then $h(G)=1$. This gives a counterexample to Shi's conjecture because $|\mu(B)|=4$, and so the minimal normal subgroup $N=Soc(G)=B^{23}$ of $G$ has $|\mu(N)|=1$. It would be interesting to generalize Mazurov's example replacing a subgroup of $Sym_{23}$ with a~subgroup of $Sym_k$, where $k$ belongs to some infinite series, possibly with another group $B$.

\begin{problem}\label{p:wreath} Does there exist a finite simple group $B$ such that for infinitely many numbers $k$, there is a finite recognizable by spectrum group $G$ with $Soc(G)=B^k$?
\end{problem}

We would like to emphasize the difference between Problem \ref{p:nonsimple}
and Problem \ref{p:wreath}: in the former, a new simple group can be chosen for every new $k$, while in the latter, we need one fixed simple group suitable for infinitely many $k$.

Observe that for all recognizable groups $G$ mentioned above,
the group $G/Soc(G)$ is solvable. We conclude our survey with
the following question posed by D. O. Revin.

\begin{problem} Does there exist a finite recognizable by spectrum group $G$ such that $G/Soc(G)$ is not solvable?
 \end{problem}

\noindent{\bf Acknowledgments.} We are grateful to A. S. Staroletov, H. He, R. Shen, A. A. Buturlakin, D. O. Revin, and S. V. Skresanov for their helpful comments and productive suggestions.

\newpage

\section*{Appendix: Simple groups with solved recognition problem}
\renewcommand{\arraystretch}{1.3}

\begin{table}[h]
\caption{\textit{Finite groups isospectral to} $L=L_{n}(q)$\\$q=p^m$, $d=(n,q-1)$, $b=((q-1)/d,m)_d$, $n_0=27$ \\$\psi,\psi_1, \chi\in\langle \varphi\rangle$, $\eta\in\langle\delta\rangle$, $|\psi|=(b)_{2'}$,  $|\psi_1|=(m)_3$, $|\chi|=(b)_2$, $|\eta|=(d)_2$}\label{tab:l}

$\begin{array}{|l|l|c|c|}
\hline
\multicolumn{2}{|c|}{\text{Conditions on }L}&h(L)&\text{Note}\\
\hline
n=2& q=9&\infty&\text{See \cite{91BrShi}}\\
& q\neq 9& 1& \\
\hline
n=3&q=3& \infty&\text{See \cite{88Shi1, 02Maz.t}}\\
&p=2 \text{ or } 3<q\equiv 3,11\Mod{12}& 1&\\
&q\equiv 5,9\Mod{12}&2& \langle \gamma \rangle\\
&q\equiv 1\Mod{6}& \tau((m)_3) &   \langle\psi\rangle\\
\hline
n=4&  p=2, 3 \text{ or } (q\equiv 1\Mod 4, m \text{ odd})&1&\\
&q\equiv 1\Mod 4, m \text{ even}, p\neq 3&2\tau((m)_2)-1&  \langle \chi\rangle\times \langle \gamma\rangle\setminus\{\gamma\}\\
& q\equiv 3,7\Mod {12}, p\neq 3&2 &   \langle\gamma\rangle\\
& q\equiv -1\Mod {12}& 2\tau((m)_3)&  \langle\psi_1\rangle \times \langle\gamma\rangle\\
\hline
n\geqslant 5, p=2& n-1=2^t&1&\\
& n-1\neq 2^t&\tau(b)&  \langle\psi\rangle\\
\hline
n\geqslant n_0 \text{ or}& n-1=p^t && \\
n\geqslant 5 \text{ prime}, &\quad n=2^s+2 &1&\\
p\equiv\kappa\Mod 4&\quad n\neq 2^s+2 &&\\
\text{for }\kappa=\pm1&\quad \quad (b)_2=1&1&\\
&\quad\quad (b)_2=2, (p-\kappa)_2\leqslant (n)_2&1&\\
&\quad\quad (b)_2=2, (p+1)_2>(n)_2&2& \langle\varphi^{m/2}\eta\rangle\\
&\quad\quad (b)_2>2\text{ or }(p-1)_2>(n)_2&2& \langle\varphi^{m/2}\gamma\eta\rangle\\
& n-1\neq p^t &&\\
& \quad (b)_2=1&&\\
\multirow{2}{*}{\framebox[2.17cm][l]{\parbox{2cm}{See \\ Lemma \ref{l:graph}}}$\ \rightarrow$}&\quad \quad\omega(L\rtimes \langle\gamma\rangle)\neq \omega(L) &\tau(b)&  \langle \psi\rangle\\
& \quad \quad\omega(L\rtimes \langle\gamma\rangle)=\omega(L) & 2\tau(b)&  \langle \psi\rangle\times \langle \gamma\rangle\\
& \quad (b)_2>1&&\\
&\quad \quad n=p^s+2^u+1&2\tau(b)-\tau((b)_{2'})&  \Theta_1 \\
&\quad \quad n\neq p^s+2^u+1&&\\
&\quad\quad\quad (n)_2\geqslant (p-\kappa)_2&3\tau(b)-3\tau((b)_{2'})&  \Theta_1\cup \Theta_2\\
&\quad\quad\quad (n)_2<(p-\kappa)_2&3\tau(b)-2\tau((b)_{2'})&  \Theta_1\cup \Theta_3^\kappa\\
\hline
\multicolumn{4}{|c|}{\Theta_1=\langle \psi\rangle\times\langle\chi\rangle\times \langle \gamma\rangle\setminus \langle \psi\rangle\gamma, \quad\Theta_2=\langle \psi\rangle\times\langle\chi^2\rangle\times \langle \gamma\eta\rangle\setminus \langle \psi\rangle\gamma\eta}\\
\multicolumn{4}{|c|}{\Theta_3^+=\langle \psi\rangle\times\langle\chi\rangle\times \langle \gamma\eta\rangle\setminus \langle \psi\rangle\gamma\eta, \quad\Theta_3^-=\langle \psi\rangle\times\langle\chi\gamma\rangle\times \langle \gamma\eta\rangle\setminus \langle \psi\rangle\gamma\eta}\\
\hline
\end{array}$
\end{table}

\begin{table}
\caption{\textit{Finite groups isospectral to} $L=U_{n}(q)$\\$q=p^m$, $d=(n,q+1)$, $b=(((q+1)/d,m)_d)_{2'}$, $n_0=27$\\$\psi,\psi_1,  \chi\in\langle \varphi\rangle$, $|\psi|=b$,  $|\psi_1|=(m)_3$, $|\chi|=2(m)_2$, $\gamma=\chi^{(m)_2}$}\label{tab:u}

$\begin{array}{|l|l|c|c|}
\hline
\multicolumn{2}{|c|}{\text{Conditions on }L}&h(L)&\text{Note}\\
\hline
n=3&q=5 \text{ or $q$ special Mersenne}& \infty&\text{See \cite{98Maz.t, 06Zav.t}}\\
&p=2 \text{ or } q\equiv 1,9\Mod{12}& 1& \\
&q\equiv 3,7\Mod{12}, q \text{ not s. M.}&2& \langle \gamma \rangle\\
&5<q\equiv 5\Mod{6}& \tau((m)_3) &  \langle\psi\rangle\\
\hline
n=4& q=2& \infty &\text{See \cite{98Maz.t}}\\
& q>2, p=2,3\text{ or } q\equiv -1\Mod 4&1& \\
& q\equiv 5,9\Mod {12}, p\neq 3&\tau(2(m)_2) &  \langle\chi\rangle\\
& q\equiv 1\Mod {12}& \tau(2(m)_6)&   \langle\psi_1\rangle\times\langle\chi\rangle\\
\hline
n=5, p=2&q=2&\infty&\text{See \cite{98Maz.t}}\\
   & q>2&1& \\
\hline
n\geqslant 6, p=2& n-1=2^t&1&\\
& n-1\neq 2^t&\tau(b)&  \langle\psi\rangle\\
\hline
n\geqslant n_0, p\neq 2& n-1=p^t &1& \\
& n-1\neq p^t &&\\
\multirow{2}{*}{\framebox[2.17cm][l]{\parbox{2cm}{See \\Lemma \ref{l:graph}}}$\ \rightarrow$}& \quad \omega(L\rtimes \langle\gamma\rangle)\neq \omega(L) &\tau(b)&   \langle \psi\rangle\\
& \quad \omega(L\rtimes \langle\gamma\rangle)=\omega(L)& &\\
& \quad \quad n\geqslant 16, (n)_2>2& 2\tau(b)&  \langle \psi\rangle\times \langle\gamma\rangle\\
& \quad \quad n\leqslant 15 \text{ or }(n)_2\leqslant 2 &\tau(2(m)_2b)&  \langle \psi\rangle\times \langle\chi\rangle\\
\hline
\end{array}$
\end{table}

\begin{table}
\caption{\textit{Finite groups isospectral to} $L=S_{2n}(q)$\\$q=p^m$, $n_0=16$,
$\chi\in\langle \varphi\rangle$, $|\chi|=(m)_2$}

$\begin{array}{|l|l|c|c|}
\hline
\multicolumn{2}{|c|}{\text{Conditions on }L}&h(L)&\text{Note}\\
\hline
n=2&q=3 \text{ or } q\neq 3^{2t+1}&\infty&\text{See \cite{98Maz.t, 00MazXuCao.t, 02Maz.t}}\\
& 3<q=3^{2t+1}&1&\\
\hline
n=3&q=2&2&O_8^+(2)\\
&q>2, p=2, 5& 1& \\
& p\neq 2,5& \tau((m)_2) &  \langle \chi\rangle\\
\hline
n=4& p\neq 7&\infty&\text{See \cite{06MazMog, 14GrSt}, \cite[Theorem 3]{18Gr.t} }\\
\hline
n\geqslant 5, p=2&\text{ none }&1& \\
\hline
n=2^s\geqslant 8\text{ or }n\geqslant n_0& 2n-1=p^t&1& \\
p\neq 2&  2n-1\neq p^t&\tau((m)_2)&  \langle \chi\rangle\\
\hline
\end{array}$

\end{table}

\begin{table}
\caption{\textit{Finite groups isospectral to} $L=O_{2n+1}(q)$\\$q=p^m$ odd, $n_0=16$, $\chi\in\langle \varphi\rangle$, $|\chi|=(m)_2$}
$\begin{array}{|l|l|c|c|}
\hline
\multicolumn{2}{|c|}{\text{Conditions on }L}&h(L)&\text{Note}\\
\hline
n=3&q=3&2&O_8^+(3)\\
&p=5& 1& \\
&q\neq 3, p\neq 5& \tau((m)_2) &   \langle \chi\rangle\\
\hline
n=4& \text{ none }&\infty& \text{See \cite{14GrSt}}\\
\hline
n=2^s\geqslant 8\text{ or }n\geqslant n_0& 2n-1=p^t&1& \\
& 2n-1\neq p^t&\tau((m)_2)&   \langle \chi\rangle\\
\hline
\end{array}$
\end{table}

\begin{table}
\caption{\textit{Finite groups isospectral to} $L=O_{2n}^+(q)$\\$q=p^m$, $n_0=19$, $\chi\in\langle \varphi\rangle$, $|\chi|=(m)_2$}\label{tab:o+}

$\begin{array}{|l|l|c|c|}
\hline

\multicolumn{2}{|c|}{\text{Conditions on }L}&h(L)&\text{Note}\\
\hline
n=4&q=2& 2&S_6(2)\\
&q=3& 2&O_7(3)\\
&q>3, p=2, 5&1&\\
&q>3, p\neq 2, 5& \tau((m)_2) &  \langle \chi\rangle\\
\hline
n\geqslant 5, p=2&\text{ none }&1&\\
\hline
n\geqslant n_0, p\neq 2&  2n-3=p^t &1& \\
&  2n-3\neq p^t &&\\
& \quad n \text{ even } &\tau((m)_2)&   \langle \chi\rangle\\
& \quad n \text{ odd, }q\equiv -1\pmod 4& 2&  \langle\gamma\rangle\\
& \quad n \text{ odd, }p\equiv 1\pmod 4& & \\
& \quad \quad 2n-3>p& \tau((m)_2)& \langle \chi\rangle\\
& \quad \quad 2n-3<p& 2\tau((m)_2)-1& \langle \chi\rangle\times \langle \gamma\rangle\setminus\{\gamma\}\\
& \quad n \text{ odd, }p\equiv -1\pmod 4, m \text{ even}& & \\
& \quad \quad 2n-3>p^2 \text{ or } 2n-3-p=2^t& \tau((m)_2)& \langle \chi\gamma\rangle\\
& \quad \quad 2n-3<p^2, 2n-3-p\neq 2^t& 2\tau((m)_2)-1&  \langle \chi\rangle\times \langle \gamma\rangle\setminus\{\gamma\}\\
\hline
\end{array}$
\end{table}

\begin{table}
\caption{\textit{Finite groups isospectral to} $L=O_{2n}^-(q)$\\$q=p^m$, $n_0=18$, $\chi\in\langle \varphi\rangle$, $|\chi|=2(m)_2$}\label{tab:o-}

$\begin{array}{|l|l|c|c|}
\hline
\multicolumn{2}{|c|}{\text{Conditions on }L}&h(L)&\text{Note}\\
\hline
n\geqslant 4, p=2&\text{ none }&1 & \\
\hline
n=2^s\geqslant 4 \text{ or }n\geqslant n_0, & 2n-3=p^t \text{ or } (4,q^n+1)=4&1& \\
p\neq 2& 2n-3\neq p^t, (4,q^n+1)=2&\tau(2(m)_2)&  \langle \chi\rangle\\
\hline
\end{array}$
\end{table}

\begin{table}[h]
\caption{\textit{Finite groups isospectral to some simple classical groups with disconnected prime graph}\\$n\geqslant 5$}\label{tab:dpg}

$\begin{array}{|l|l|c|c|}
\hline
L&\multicolumn{1}{c|}{\text{Conditions on }L}&h(L)&\multicolumn{1}{c|}{\text{Note}}\\
\hline
U_n(3)& n\text{ is a prime}&2&\langle \gamma\rangle\\
& n-1\text{ is a prime}, (n)_2>4&2&\langle \gamma\rangle \\
& n-1\text{ is a prime}, (n)_2=4&1& \\
\hline
S_{2n}(3), O_{2n+1}(3)&n\text{ is a prime}&1& \\
\hline
O_{2n}^+(3)&n\text{ is a prime}&2&\langle\gamma\rangle \\
&n-1\text{ is a prime}&1&\\
\hline
O_{2n}^+(5)&n\text{ is a prime}&1& \\
\hline
O_{2n}^-(3)&n\text{ is a prime or }n=2^m+1&1& \\
\hline
\end{array}$
\end{table}

\begin{table}[h]
\caption{\textit{Finite groups isospectral to exceptional groups of Lie type}\\ $q=p^m$}\label{tab:ex}

$\begin{array}{|m{2.1cm}|l|c|c|}
\hline
$L$&\multicolumn{1}{c|}{\text{Conditions on }L}&h(L)&\multicolumn{1}{c|}{\text{Note}}\\
\hline
$^3D_4(q)$&q=2&\infty&\text{See \cite{13Maz.t}}\\
&q>2, p\in\{2,3,7,11\}&1& \\
&p\not\in\{2,3,7,11\}& \tau((m)_2)& \langle \varphi^{(m)_{2'}}\rangle\\
\hline
$F_4(q)$&p\in\{2,3,7,11\}&1& \\
&p\not\in\{2,3,7,11\}& \tau((m)_2)& \langle \varphi^{(m)_{2'}}\rangle\\
\hline
$E^\varepsilon_6(q)$&p\in\{2,11\} \text{ or }(3,q-\varepsilon)=1&1& \\
&p\not\in\{2,11\}, (3,q-\varepsilon)=3& \tau((m)_3)& \langle \varphi^{(m)_{3'}}\rangle\\
\hline
$E_7(q)$&p\in\{2,13,17\}&1& \\
&p\not\in\{2,13,17\}& \tau((m)_2)&  \langle \varphi^{(m)_{2'}}\rangle\\
\hline
Other \hspace{7mm} exceptional groups&\text{none }&1& \\
\hline
\end{array}$
\end{table}

\begin{table}[h]
\caption{\textit{Finite groups isospectral to alternating and sporadic groups}}\label{tab:alt}

$\begin{array}{|l|l|c|c|}
\hline
L&\text{Conditions on }L&h(L)&\text{Note}\\
\hline
Alt_n&n\neq 6,10&1&\\
&n=6&\infty& \text{See \cite{91BrShi} }\\
&n=10&\infty& \text{See \cite{98Maz.t}}\\
\hline
\text{Sporadic group}&L\neq J_2&1&\\
&L=J_2&\infty&\text{See \cite{94PrShi, 98MazShi}}\\
\hline
\end{array}$
\end{table}

\begin{table}[p]
\caption{\emph{Unrecognizable simple groups $L$ and finite groups
$G$ such that $\omega(G)=\omega(L)$ and $G$ is not an almost simple group with socle $L$}}\label{tab:un}

\begin{tabular}{|l|l|c|c|c|c|l|}
\hline
\multicolumn{2}{|l|}{$L$}& \multicolumn{3}{c|}{$G$}& Examples&References\\
\cline{3-5}
\multicolumn{2}{|l|}{}& $\pi(K)$&$S$&$G/H$& &\\
\hline
\multicolumn{2}{|l|}{$Alt_6$}& $\{2\}$&$Alt_5$&$1$& ${2^4:Alt_5}$&\cite{91BrShi, 13LytY}\\
\hline
\multicolumn{2}{|l|}{$Alt_{10}$}&$\{2,3,7\}$&$Alt_5$&$2$&$(7^4\times 3^{12}):(2.Sym_5)$&\cite{98Maz.t, 10Sta.t, 15LytY.t}\\
\hline
\multicolumn{2}{|l|}{$J_2$}&$\varnothing$&$Alt_8$&$2$&$Sym_8$&\cite{94PrShi, 98MazShi, 15LytY.t}\\
\cline{3-6}
\multicolumn{2}{|l|}{ }&$\{2\}$&$Alt_8$&$1$&$2^6:Alt_8$&\\
\hline
\multicolumn{2}{|l|}{$^3D_4(2)$}&$\{2\}$&$^3D_4(2)$&1&$2^{24}:{}^3D_4(2)$&\cite{13Maz.t}\\
\hline
\multicolumn{2}{|l|}{$L_3(3)$}&\multicolumn{3}{c|}{solvable Frobenius group}&$13^4:(2.Sym_4)$&\cite{88Shi1, 02Maz.t, 13LytY}\\
\hline
\multicolumn{2}{|l|}{$U_3(3)$}&$\{2\}$&$U_3(3)$& $1,2$& $2^6:U_3(3)$&\cite{98Maz.t, 03Ale.t, 15Maz1.t,17Lyt}\\
\cline{3-6}
\multicolumn{2}{|l|}{ }&$\{2\}$&$L_2(7)$&1& $4^3.L_2(7)$ &\\
\cline{3-6}
\multicolumn{2}{|l|}{ }&$\{2\}$&$L_2(7)$&2& $2^6:PGL_2(7)$&\\
\cline{3-6}
\multicolumn{2}{|l|}{ }&\multicolumn{3}{c|}{2-Frobenius group}&$2^{18}:(7:3)$&\\
\cline{3-6}
\multicolumn{2}{|l|}{ }&\multicolumn{3}{c|}{solvable Frobenius group}&$7^4:(3:8)$&\\
\hline
\multicolumn{2}{|l|}{$U_3(5)$}&$\subseteq\{2\}$&$L_3(4)$&$1,2$&$2^9:L_3(4)$, $L_3(4).\langle \gamma\rangle$ &\cite{98Maz.t, 02Ale.t}\\
\cline{3-6}
\multicolumn{2}{|l|}{ }&$\{2\}$&$Alt_7$&$1$&unknown&\\
\hline
\multicolumn{2}{|l|}{$U_3(q)$, $q\geqslant 7$,}& $\{2\}$&$U_3(q)$&$1,2$&$2^{q^2-q}:U_3(q)$&\cite{98Maz.t, 06Zav.t, 02Ale.t}\\
\multicolumn{2}{|l|}{$q$ special Mersenne}&&&&&\\
\hline
\multicolumn{2}{|l|}{$U_5(2)$}&$\{3\}$&$U_5(2)$&$1$&$3^{10}:U_5(2)$&\cite{98Maz.t, 11Gr}\\
\cline{3-6}
\multicolumn{2}{|l|}{ }&$\{3\}$ &$M_{11}$&$1$&$3^5:M_{11}$  &\\
\hline
\multicolumn{2}{|l|}{$S_4(3)$}
&$\{2,3\}$&$Alt_5$&$1$&$(2^4\times 3^4):Alt_5$&\cite{02Maz.t, 10Zav.t, 19Lyt}\\
\cline{3-6}
\multicolumn{2}{|l|}{ }&$\{3\}$&$Alt_5$&$2$&$3^4:Sym_5$&\\
\cline{3-6}
\multicolumn{2}{|l|}{ }&\multicolumn{3}{c|}{2-Frobenius group}&$[3^{24}]:(5:4)$&\\
\hline
$S_4(q)$,  &$p=2$&$\subseteq\{2\}$ &$L_2(q^2)$&$\leq\langle \alpha^m\rangle$&$2^{8m}:L_2(q^2)$, &\cite{00MazXuCao.t, 18Lyt}\\
$q=p^m$ &&&&&$L_2(q^2).\langle\alpha^m\rangle$&\\
\cline{2-7}
&$p=3$, &
 $\{3\}$ &$L_2(q^2)$&$\leq\langle \alpha\rangle$&$3^{28m}:L_2(q^2)$&\cite{02Maz.t, 18Lyt}\\
&$m$ even&&&&&\\
\cline{2-7}
&$p\geqslant 5$&$\{p\}$&$L_2(q^2)$&$\leq\langle \alpha\rangle$&$p^{8m}:(L_2(q^2).\langle \alpha^m\rangle)$&\cite{02Maz.t, 18Lyt}\\

\hline
$S_8(q)$,&$p\neq 2,7$&$\{p\}$&$O_8^-(q)$&$\leq\langle\chi\rangle$&$p^{8m}:(O^-_8(q).\langle\gamma\rangle)$&\cite{18Gr.t}\\
$q=p^m$&&&&&&\\
\hline
$O_9(q)$, &$p=2$&$\subseteq \{2\}$&$O_8^-(q)$&$2$-group&$2^{8m}:O^-_8(q)$, $O_8^-(q).\langle \gamma\rangle$& \cite{06MazMog, 14GrSt,  15VasGr1}\\
\cline{2-7}
$q=p^m$&$p\neq 2$&$\{p\}$&$O_8^-(q)$&$\leq\langle\chi\rangle$&$p^{8m}:O^-_8(q)$&\cite{09VasGorGr.t, 14GrSt, 18Gr.t}\\
\hline
\multicolumn{7}{|c|}{$\alpha$ is a field automorphism of $L_2(q^2)$, $q=p^m$, of order $2(m)_2$}\\
\multicolumn{7}{|c|}{$\gamma$ is the graph automorphism of $S$ as defined in Subsection \ref{ss:def} }\\
\multicolumn{7}{|c|}{$\chi$ is the automorphism of $O_8^-(q)$ as in Table \ref{tab:o-}; in particular, $\chi^m=\gamma$}\\
\hline
\end{tabular}

\end{table}

\end{document}